\documentclass{amsart}
\usepackage{amsthm,amsmath,amsfonts}
\usepackage{mathrsfs}
\usepackage[colorlinks=false,linktocpage=true]{hyperref}
\usepackage{tocvsec2}
\usepackage{latexsym}
\numberwithin{equation}{section}
\font\tenscrpt=eusm10
\font\sevenscrpt=eusm10 scaled 700
\font\fivescrpt=eusm10 scaled 500
\newfam\eusmfam
\textfont\eusmfam=\tenscrpt
\scriptfont\eusmfam=\sevenscrpt
\scriptscriptfont\eusmfam=\fivescrpt
\def\scr#1{{\fam\eusmfam\relax#1}}
\font\eight=cmr8

\def\qed{\quad\vcenter{\hrule\hbox{\vrule height.6em\kern.6em\vrule}\hrule}}
\newenvironment{pf}{{\textsc Proof.\quad}}{$\qed$\bigskip\newline}
\newenvironment{pf*}[1]{{\textsc #1.\quad}}{$\qed$\bigskip\newline}
\newtheorem{thm}{Theorem}[section]

\newtheorem{lem}{Lemma}[section]

\newcommand{\thmref}[1]{Theorem~{\rm $\ref{#1}$}}

\def\eqdef{\overset{\triangle}{=}}

\def\eqdef{\overset{\triangle}{=}}
\def\utx{u(t,x)}
\def\vsx{v(s,x)}

\def\Xx{X^x}
\def\Xix{X^{-ix}}

\def\KSfBTP{{\mathbb {A}^{f,X}_{B}(t,x)}}
\def\KSfBTPz{{\mathbb {A}^{f,X}_{B}(0,x)}}

\def\P{{\mathbb P}}
\def\GE{{\mathbb {E}}^{\mathbb {C}}}
\def\EP{{\mathbb E}_{\P}}

\def\R{{\mathbb R}}
\def\Rd{{\mathbb R}^d}

\def\Rp{{\mathbb R}_+}
\def\Rpmo{{\overset{\circ}{{\mathbb R}}_{\pm}}}
\def\Rpo{{\overset{\circ}{{\mathbb R}}_{+}}}
\def\Rmo{{\overset{\circ}{{\mathbb R}}_{-}}}

\def\P{{\mathbb P}}

\def\EP{{\mathbb E}_{\P}}

\def\R{{\mathbb R}}

\begin{document}
\title[LKS via imaginary-Brownian-time-Brownian-angle processes]{A linearized Kuramoto-Sivashinsky PDE via an imaginary-Brownian-time-Brownian-angle process}
\author{Hassan Allouba}
\address{Department of Mathematical Sciences, Kent State University, Kent, Ohio 44242}
\email{allouba@math.kent.edu}
\date{5/20/2002}
 \maketitle
\begin{abstract}
We introduce a new imaginary-Brownian-time-Brownian-angle process, which we also call the linear-Kuramoto-Sivashinsky process (LKSP).
Building on our techniques in two recent articles involving the connection of Brownian-time processes to fourth order PDEs, we give an explicit
solution to a linearized Kuramoto-Sivashinsky PDE in $\scriptstyle d$-dimensional space: $\scriptstyle{\eight{ u_t=-\frac18\Delta^2u-\frac12\Delta u-\frac12u}}$.
 The solution is given in terms of a functional of our LKSP.
 \end{abstract}

 \section{Statements and discussions of results.}
One of the prominent equations in modern applied mathematics is the celebrated Kuramoto-Sivashinsky (KS) PDE.
This nonlinear equation has generated  a lot of interest in the PDE literature (see e.g., \cite{CT,CF,FK,FT,TW} and many other papers).  In the field of
stochastic processes, a
great deal of interest is directed at the study of processes in which time is replaced in one way or another by a Brownian motion, and this interest has picked up
considerably (see e.g., \cite{AZ,AFK,BeRoVa,Bu1,Bu2,ZS,ZSY,KL,KL2,F,HO}) after the fundamental work of Burdzy on iterated Brownian motion (\cite{Bu1,Bu2}).  In
\cite{AZ,AFK}, we provided a unified framework for such iterated processes (including the IBM of Burdzy) 
and introduced several interesting new ones,
through a large class of processes that we called Brownian-time processes.  We then related them to different fourth order PDEs.   In this article,
and as announced in \cite{AFK}, we modify our process in Theorem 1.2 \cite{AFK} and build on our methods in \cite{AFK} to give an explicit
solution to a linear version of the KS PDE.  One modification needed is the introduction of $i=\sqrt{-1}$ in both the Brownian-time and the
Brownian-exponential, and that leads to a new process we call imaginary-Brownian-time-Brownian-angle process IBTBAP, starting at
$f:\Rd\to\R$:
 \begin{equation}
\KSfBTP\eqdef\begin{cases} f(\Xx(iB(t))) \exp\left(iB(t)\right),  & B(t)\ge0;\cr
f(i\Xix(-iB(t)))\exp\left(iB(t)\right),  & B(t)<0;
\end{cases}
\label{KSprocess}
\end{equation}
where $\Xx$ is an $\Rd$-valued Brownian motion starting from $x\in\Rd$, $\Xix$ is an independent $i\Rd$-valued BM starting at $-ix$ (so that $i\Xix$ starts
at $x$), and both are independent of the inner standard $\R$-valued Brownian motion $B$ starting from $0$.    The time of the outer Brownian motions $\Xx$ and $\Xix$ is replaced
by an imaginary positive Brownian time; and, when $f$ is real-valued as we will assume here, the angle of $\KSfBTP$ is the Brownian motion $B$.
We think of the imaginary-time processes $\{\Xx(is),s\ge0\}$ and $\{i\Xix(-is),s\le0\}$ as having the same complex Gaussian distribution on $\Rd$ with the corresponding complex distributional density
$$p_{is}^{(d)}(x,y)=\frac{1}{({2\pi i s})^{d/2}}e^{-|x-y|^2/2is}.$$
We will also call the process given by \eqref{KSprocess} the $d$-dimensional
Linear-Kuramoto-Sivashinsky process (LKSP) starting at $f$ (clearly $\KSfBTPz=f(x)$).    The dimension in $d$-dimensional IBTBAP (or $d$-dimensional LKSP) refers to the dimension of the BMs $\Xx$ and $\Xix$, which is also the
dimension of the spatial variable in the associated linearized KS PDE as we will see shortly.

Now, motivated by the definitions of $v_\epsilon$ and
$u_\epsilon$ in the proof of Theorem 1.2 in \cite{AFK}, we let
\begin{equation}
\begin{split}
&v(s,x)\eqdef\exp\left(is\right)\int_{\Rd}f(y)\frac{1}{{(2\pi is)}^{d/2}}e^{-|x-y|^2/2is} dy\\
&u(t,x)\eqdef \int_{-\infty}^0v(s,x)p_t(0,s)ds+ \int_0^\infty v(s,x) p_t(0,s)ds
\end{split}
\label{expectn}
\end{equation}
where $p_t(0,s)$ is the transition density of the inner (one-dimensional) Brownian motion $B$:
$$p_t(0,s)=\frac{1}{\sqrt{2\pi t}}e^{-s^2/2t}.$$
We may think of $v$ and $u$ in terms of complex expectation by defining
$v(s,x)\eqdef\GE\left[f(\Xx(is))\exp\left(is\right)\right]$ and $u(t,x)\eqdef\GE\left[\KSfBTP\right]$.
A more detailed study of the rich connection between our process and its complex distribution to the KS PDE and its implications is the subject of an upcoming
article \cite{A3}.
We are now ready to state our main result.
\begin{thm}
Let $f\in C^2_c(\Rd;\R)$ with $D_{ij} f$ H\"older continuous with exponent $0<\alpha\le1$, for all $1\le i,j\le d$.
If $\utx$ is given by \eqref{expectn} then $\utx$ solves the linearized
Kuramoto-Sivashinsky PDE
\begin{equation}
\begin{cases}
 \dfrac{\partial}{\partial t} \utx=- \dfrac{1}{8}\Delta^2\utx-\dfrac12\Delta \utx-\dfrac{1}{2}\utx,
 &t>0,\,x\in\Rd;
\cr u(0,x)=f(x), &x\in\Rd.
\end{cases}
\label{KSPDE}
\end{equation}
\label{PDEKS}
\end{thm}

\section{Proof of the main result}
\begin{pf*}{Proof of \thmref{PDEKS}}
Let $u$ and $v$ be as given in \eqref{expectn}.   Differentiating $\utx$ with respect to $t$ and putting the derivative under the integral,
which is easily justified by the dominated convergence theorem, then using the fact that $p_t(0,s)$ satisfies the heat equation
$$\frac{\partial}{\partial t} p_t(0,s)=\frac12 \frac{\partial^2}{\partial s^2}p_t(0,s)$$
and integrating by parts twice using the fact that the boundary terms vanish at $\pm\infty$ and that $(\partial/\partial s)p_t(0,s)=0$ at $s=0$,  we obtain
\begin{equation}
\begin{split}
\frac{\partial}{\partial t} \utx & = \int_{-\infty}^0 \vsx\frac{\partial}{\partial t}p_t(0,s) ds+\int_0^\infty \vsx\frac{\partial}{\partial t}p_t(0,s) ds
\\&= \frac12\left[\int_{-\infty}^0\vsx\frac{\partial^2}{\partial s^2} p_t(0,s) ds+\int_0^\infty\vsx\frac{\partial^2}{\partial s^2} p_t(0,s) ds\right]\\
&=\frac12p_t (0,0)\left[\left.\left( {\frac {\partial }{\partial s}}\vsx\right)\right|_{ s=0^-} +\left.\left( {\frac {\partial }{\partial s}}\vsx\right)\right|_{ s=0^+}\right]
\\&+\frac12\int_{-\infty}^0 p_t(0,s)\frac{\partial^2}{\partial s^2} \vsx  ds+\frac12\int_0^\infty p_t(0,s)\frac{\partial^2}{\partial s^2} \vsx  ds\\
&=\frac12\int_{-\infty}^0 p_t(0,s)\left[-\frac{1}{4}\Delta^2v \left( s,x \right) -
\Delta v \left( s,x \right) -v \left( s,x \right)\right] ds
\\&+\frac12\int_0^\infty p_t(0,s)\left[-\frac{1}{4}\Delta^2v \left( s,x \right) -
\Delta v \left( s,x \right) -v \left( s,x \right)\right] ds
\\ &=-\frac18\Delta^2\utx-\frac12\Delta \utx- \frac12\utx
\end{split}
\label{tder77}
\end{equation}
 where for the last two equalities in \eqref{tder77} we have used the fact that
\begin{equation}
\begin{split}
\frac{\partial v}{\partial s}&=\frac{i}{2}\Delta v \left( s,x \right) +iv \left( s,x \right)\\
\frac{\partial^2 v}{\partial s^2}&=-\frac{1}{4}{\Delta^2}v \left( s,x \right) -
\Delta v \left( s,x \right) -v \left( s,x \right),
\end{split}
\label{vint}
\end{equation}
and the conditions on $f$  to take the applications of the derivatives outside the integrals in \eqref{tder77} and \eqref{vint} (the steps of Lemma 2.1 in \cite{AFK}
easily translates to our setting here,  see the discussion below).
Clearly $u(0,x)=f(x)$, and the proof is complete.
\end{pf*}
As we indicated above, only minor changes to Lemma 2.1 in \cite{AFK} are needed to justify pulling the derivatives outside the integrals in
\eqref{tder77} under the conditions on $f$ of  \thmref{PDEKS}.   We now adapt Lemma 2.1 \cite{AFK} to our setting here, and we point out
the necessary changes in its proof:
\begin{lem} Let $v(s,x)$ be given by \eqref{expectn} and let  $f$ be as in \thmref{PDEKS}.  Let
\begin{equation}
\begin{split}
u_1(t,x)\eqdef\int_{-\infty}^0 v(s,x) p_t(0,s) ds \quad\mbox{and}\quad u_2(t,x)\eqdef\int_0^\infty v(s,x) p_t(0,s) ds,
\end{split}
\label{diffunderint}
\end{equation}
then  $\Delta^2u_1(t,x)$ and $\Delta^2u_2(t,x)$ are finite and
\begin{equation}
\begin{split}
\Delta^2u_1(t,x)=\int_{-\infty}^0\Delta^2 v(s,x) p_t(0,s) ds \quad\mbox{and}\quad \Delta^2u_2(t,x)=\int_0^\infty\Delta^2 v(s,x) p_t(0,s) ds.
\end{split}
\label{diffunderint1}
\end{equation}
\end{lem}
\begin{pf}
As in the proof of Lemma 2.1 \cite{AFK}, letting $\Rpo=(0,\infty)$ and $\Rmo=(-\infty,0)$, it suffices to show
\begin{equation}
\frac{\partial^4 }{\partial x_j^4}\int_\Rpmo v(s,x) p_t(0,s) ds=\int_\Rpmo\frac{\partial^4}{\partial x_j^4}v(s,x) p_t(0,s) ds, \qquad j=1,\ldots, d.
\label{eq1}
\end{equation}
Letting $p_{is}^{(d)}(x,y)=(2\pi is)^{-d/2}e^{-|x-y|^{2}/2is}$ and using the conditions on $f$, we easily get
\begin{equation}
\begin{split}   \label{fstep}
\frac{\partial^4}{\partial x_j^4}v(s,x) p_t(0,s)&= \exp\left(is\right)\left(\int_{\Rd}f(y)\frac{\partial^4}{\partial y_j^4}p^{(d)}_{is}(x,y)dy\right)p_t(0,s)\\&=\exp\left(is\right)
\left(\int_{\Rd}\frac{\partial^2}{\partial y_j^2}f(y)\frac{\partial^2}{\partial y_j^2}p^{(d)}_{is}(x,y)dy\right)p_t(0,s). \end{split}
\end{equation}
Rewriting the last term in \eqref{fstep}, and letting $h_j(y)\eqdef\partial^2 f(y)/\partial y_j^2$, we have
\begin{equation}
\begin{split}
&\left|\exp\left(is\right)\left(\int_{\Rd}(2\pi is)^{-{d/2}}\left(\frac{-(x_j-y_j)^2+is}{s^2}\right)e^{-|x-y|^2/2is}h_j(y)dy\right)\frac{e^{-s^2/2t}}{\sqrt{2\pi t}}\right|\\
&= \frac{e^{-s^2/2t}}{\sqrt{2\pi t}}\left|\left(\int_{\Rd}(2\pi is)^{-{d/2}}\left(\frac{-(x_j-y_j)^2+is}{s^2}\right)e^{-|x-y|^2/2is}(h_j(y)-h_j(x))dy\right)\right|\\
&\le\frac{e^{-s^2/2t}}{\sqrt{2\pi t}}\int_{\Rd}(2\pi|s|)^{-{d/2}}\left|\frac{-(\tilde{x}_j-y_j)^2+|s|}{s^2}\right|e^{-|\tilde{x}-y|^2/2|s|}\left|h_j(y)-h_j(\tilde{x})\right|dy\\
&=\frac{e^{-s^2/2t}}{\sqrt{2\pi t}}\EP\left|\left(\frac{(\tilde{x}_j-W_j^{\tilde{x}}(|s|))^2-|s|}{s^2}\right)\left(h_j(W^{\tilde{x}}(|s|))-h_j(\tilde{x})\right)\right|,
\end{split}
\label{sstep}
\end{equation}
for some $\tilde{x}\in\Rd$ where $\tilde{x}_j=\pm x_j$ for $j=1\ldots, d$; and where $W^{\tilde{x}}:\Omega\times\Rp\to\Rd$ is a standard Brownian motion
starting at $\tilde{x}\in\Rd$ on a probability space $(\Omega, \scr{F}, \P)$, and
$W_j^{\tilde{x}}$ is its $j$-th component.    The inequality in \eqref{sstep} follows easily if  $h_j$ is a polynomial, and standard approximation
yields the inequality for $h_j\in C_c(\Rd;\R)$.  Now, exactly as in \cite{AFK} $(2.6)$ and $(2.7)$; we use the Brownian motion scaling for $W^{\tilde{x}}$, the
Cauchy-Shwarz inequality on the last term in \eqref{sstep}, and the H\"older condition on $h_j$ to deduce that the last term in \eqref{sstep} is
bounded above by
$K\exp{(-s^2/2t)}/(\sqrt{2\pi t}|s|^{1-\alpha/2})\in L^1((-\infty,0),ds)\cap\, L^1((0,\infty),ds)$; hence
$\left|{\partial^4}/{\partial x_j^4}v(s,x) p_t(0,s)\right|\in L^1((-\infty,0),ds)\cap L^1((0,\infty),ds)$, which completes the proof by standard analysis.
\end{pf}
{\bf Acknowledgements.}  I'd like to thank Ciprian Foias for his encouragement to pursue this project and for his support.  I also enjoyed several
one on one fruitful discussions with him.  This research is supported in part by NSA grant MDA904-02-1-0083.

\end{document}